\documentclass[11pt,
]{article}

\usepackage{amssymb,amstext,amsthm,latexsym}
\usepackage{amsmath}
\usepackage{amsfonts}
\usepackage{mathrsfs}
\usepackage{mparhack}
\usepackage[bbgreekl]{mathbbol}

\input{cha.sty}

\begin{document}

\title
%[Countable cofinality of definable chains in partial orders] 
{On countable cofinality of definable chains in 
Borel partial orders
\thanks{Partial support of ESI 2013 Set theory program and 
RFFI grant 13-01-000006 acknowledged.}}

\author 
%, Lyubetsky]
{Vladimir~Kanovei\thanks{IITP RAS and MIIT,
  Moscow, Russia, \ {\tt kanovei@googlemail.com}}  
%\and
%  Vassily~Lyubetsky\inst{1} 
}

\date 
{\today}

\maketitle

\begin{abstract}
We prove that in some cases definable chains of Borel
partial orderings are necessarily countably cofinal.
This includes the following cases:
analytic chains,
\ROD\ chains in the Solovay model,
and $\fs12$ chains in the assumption that
$\omi^{\rL[x]}<\omi$ for all reals $x$.
\end{abstract}

\subsection*{Introduction}
 
Studies of maximal chains in partially ordered sets go back 
to as early as Hausdorff \cite{h07,h09}, where this issue 
appeared in connection with Du Bois Reymond's investigations 
\cite{dbr',dbR} of orders of infinity.
Using axiom of choice, Hausdorff proved the existence of
maximal chains (called \rit{pantachies}) in any partial
ordering.
On the other hand, Hausdorff clearly understood the difference
between such a pure existence proof and an actual construction
of a maximal chain --- see \eg\ \cite[p.~110]{h07} or
comments in \cite{fi} --- which we would call now the existence
of \rit{definable} maximal chains.

The following theorem is the main content of this note.
It shows that in some notable cases definable chains
are necessarily countably cofinal.

\bte
\lam{mt}
If\/ $\le$ is a Borel PQO on a (Borel) set\/ $D=\dom{(\le)}$,  
$X\sq D$, and\/ ${\le}\res X$ is a linear quasi-order
{\rm(= chain)},  
then\/ $\ang{X;{\le}}$ is countably cofinal in each of the 
following three cases: 
\ben
\renu
%\tenu{\arabic{enumi}}
\itla{mt1}
$X$ is a\/ {\  $\fs11$} set --- 
and in this case, moreover, 
there is no strictly increasing\/ $\omi$-sequences 
in\/ $X$,   

\itla{mt2}
$X$ is a\/ \ROD\ set in the Solovay model, 
%{\rm(= Theorem~\ref{teb})},

\itla{mt3}
$X$ is a\/ $\fs12$ set, and\/
$\omi^{\rL[r]}<\omi$ for every real\/ $r$. 
%\ {\rm(= Theorem~\ref{tec})}.
\een
Therefore, if, in addition, it is known that\/ $\stk D\le$
does not have maximal chains of countable cofinality, then
in all three cases\/ $X$ is not a maximal chain.
\ete 

Part \ref{mt1} is proved by reduction to a result
in \cite{hms}.
Part \ref{mt2} is already known from \cite{kl}, but we
present here a simplified proof in order to make the
exposition self-contained, since the result is used in
the proof of \ref{mt3}.

The additional condition in the theorem,
of uncountable cofinality of all maximal chains,
holds for many partial orders of interest, \eg,
the eventual domination order on sets like $\om^\om$ or
$\dR^\om,$ or the \rit{rate of growth order} defined
on $\dR^\om$ by
$$
x<_{\text{\sc rg}}y
\quad\text{iff}\quad
\lim_{n\to\infty}\frac{x(n)}{y(n)}=\infty\,.
$$
Needless to say that chains, gaps, and similar structures
related to these or similar orderings have been subject of
extended studies, of which we mention \cite{fa,fa1,to,hom}
among those in which the definability aspect is considered.

We end the introduction with a review of basic notation
related to orderings.
\bde
\item[\rm PQO, \it partial quasi-order$:$] \ 
$x\le x$ and $x\le y\land y\le z\imp x\le z$  
in the domain; 

\item[\rm LQO,  
\it linear quasi-order\,$:$] \  
in addition, $x\le y\lor y\le x$  
in the domain; 

\item[\rm LO, \it linear order\,$:$]  \ 
in addition, $x\le y\land y\le x\imp x=y$  
in the domain; 

\item[\it sub-order\,$:$] \ 
restriction of the given PQO to a subset of
its domain.

\item[$\lexe$$:$] \ 
the {\it lexicographical\/} order 
on sets of the form $2^\xi,\;\xi\in\Ord$.
\ede

\subsection{Analytic linear suborders of Borel PQOs}
\label{sek1}

In this Section, we prove Theorem \ref{mt}\ref{mt1}. 
Thus suppose that $\le$ is a Borel PQO on a Borel set 
$D\sq\bn$,  
$X\sq D$ is a $\fs11$ set, and\/ ${\le}\res X$ is a 
{\it  linear\/} quasi-order.  
Prove that this restricted 
quasi-order $\ang{X;{\le}}$ has no strictly increasing 
\ddd\omi chains. 

The proof is based on the following well-known lemma.

\ble 
\lam{lea}
Every Borel LQO\/ $\le$ is countably cofinal, 
and moreover, there is no strictly increasing\/ $\omi$-sequences.
\ele 
\bpf[lemma]
By a result in Harrington -- Marker -- Shelah \cite{hms}, 
there is an ordinal $\xi<\omi$ and a Borel map 
$f:X=\dom{(\le)}\to 2^\xi$ such that 
$$
x \le y
\quad\text{iff}\quad
f(x) \lexe f(y)
$$ 
for all $x, y \in X$. 
But the lemma easily holds for 
$\ang{2^\xi;\lexe}$.
\epf

Coming back to the proof of Theorem \ref{mt}\ref{mt1}, we 
describe the idea: 
find a Borel set $W\sq D$ such that $X\sq W$ and still 
${\le}\res W$ is linear, then use Lemma \ref{lea}.
We find such a set by means of the following two-step 
procedure.

Note that the set $Y$ of all elements in $D$, 
$\le$-comparable with every element $x \in X$, 
is $\fp11$,  
and $X \sq Y$ (as $\le$ is linear on $X$).
By the Luzin Separation theorem, there is a Borel set 
$Z$ such that $X \sq Z \sq Y$. 
This ends step 1. 

Now, at the 2nd step,  the set $U$ of all elements in 
$Z$, comparable with every element in $Z$, is $\fp11$, 
and we have $X \sq U$. 
Once again, by Separation, there is a Borel set $W$ such 
that $X \sq W \sq U$.

By construction, $\le$ is linear on $U$, and hence on $W$. 
Therefore, there is no increasing \ddd\omi sequence in $W$ 
by Lemma~\ref{lea}. 
But $X\sq W$.

\qeDD{Theorem \ref{mt}\ref{mt1}}\vtm

The next immediate corollary says that maximal chains cannot 
be analytic provided they cannot be countably cofinal. 

\bcor
\lam{cora}
If\/ $\le$ is a Borel PQO, 
and 
every countable set\/ $D'\sq \dom{\le}$ has a strict upper bound,
then there is no maximal\/ $\fs11$ chains 
in\/ $\ang{D;\le}$.\qed
\ecor

\subsection{Near-counterexamples}
\label{sec2}

The following examples show that  
Theorem \ref{mt}\ref{mt1} is not true any more 
for different extensions of the domain of 
$\fs11$ suborders of a Borel partial quasi-orders, 
such as $\fs11$ and $\fp11$ linear quasi-orders --- 
not necessarily suborders of Borel orderings, as well as 
$\fd12$ and $\fp11$ suborders of Borel orderings.
In each of these classes, a counterexample of cofinality 
$\omi$ will be defined.

\bex
[$\fs11$ LQO] 
\lam{ex1}
Consider a recursive coding of sets of rationals by reals. 
Let $Q_x$ be the set coded by a real $x$.
Let $X_\al$ be the set of all reals $x$ such that the 
maximal well-ordered initial segment
of $Q_x$ has the order type $\al$. 
We define 
$$
x\le y
\quad\text{iff}\quad
\sus\al\,\sus\ba\:(x\in X_\al\land y\in X_\ba\land\al\le\ba).
$$
Then $\le$ is a $\is11$ LQO of cofinality $\om_1$.\qed
\eex

\bex
[$\fp11$ LQO] 
\lam{ex2}
Let $D\sq\bn$ be the $\ip11$ set of codes of (countable) ordinals. 
Then
$$
x\le y
\quad\text{iff}\quad
x,y\in D\,\land\, |x|\le |y|
$$ 
is a $\ip11$ LQO of cofinality $\omi$.
\qed
\eex

\bex
[$\fp11$ LO] 
\lam{ex3}
To sharpen Example \ref{ex2}, define
$$
x\le y
\quad\text{iff}\quad
x,y\in D\;\land\;
\big( {|x|<|y|}\,\lor\, 
{(|x|=|y|\land x\lex y)}\big);
$$ 
this is a $\ip11$ LO of cofinality $\omi$.\qed
\eex

\bex
[$\fd12$ suborders] 
\lam{ex4}
Let $\le$ be the eventual domination order on $\om^\om$. 
Assuming the axiom of constructibility $\rV=\rL$, one 
can define a 
strictly $\le$-increasing $\id12$ $\omi$-sequence 
$\sis{x_\al}{\al<\omi}$ in $\om^\om$.\qed 
\eex

\bex
[$\fp11$ suborders] 
\lam{ex5}
Define a PQO $\le$ on $(\om\bez\ans0)^\om$ so that
$$
x\le y
\quad\text{iff}\quad
\lim_{n\to\iy}\:y(n)\,{/\hspace{-1ex}/}\,x(n)=\iy\,.
$$
Assuming the axiom of constructibility $\rV=\rL$, define a 
strictly $\le$-increasing $\id12$ $\omi$-sequence 
$\sis{x_\al}{\al<\omi}$ in $\om^\om$. 
By the Novikov -- Kondo -- Addison $\ip11$ Uniformization 
theorem, there is a $\ip11$ set 
$\sis{\ang{x_\al,y_\al}}{\al<\omi}$, 
such that $\al\ne\ba\,\imp\,y_\al\ne y_\ba$, 
and we may assume that each $y_\al$ belongs to $2^\om$.

Let $z_\al(n)=2^{x_\al(n)}\cdot 3^{y_\al(n)}$, $\kaz n$.
Then the \ddd\omi sequence $\sis{z_\al}{\al<\omi}$ is $\ip11$ 
and strictly $\le$-increasing.\qed
\eex

\subsection{Definable linear suborders 
in the Solovay model}

Here we prove Theorem \ref{mt}\ref{mt2}. 
Arguing in the Solovay model 
(a model of \ZFC\ defined in \cite{sol}, in which all \ROD\ 
sets of reals are Lebesgue measurable), 
we assume that $\le$ is a Borel PQO on a Borel set 
$D\sq\bn$,  
$X\sq D$ is a $\ROD$ (real-ordinal definable) set, 
and\/ ${\le}\res X$ is a 
{\it  linear\/} quasi-order.  

Prove that the restricted quasi-order $\ang{X;{\le}}$ is 
countably cofinal. 

It is known that in the Solovay model any \ROD\ set in $\bn$ 
is a union of a \ROD\ \ddd\omi sequence of analytic sets. 
Thus there is a \ddd\sq increasing \ROD\ sequence 
$\sis{X_\al}{\al<\omi}$ of $\fs11$ sets $X_\al$, such that 
$X=\bigcup_{\al<\omi}X_\al$. 

As the sets $X_\al$ are ctbly $\le$-cofinal by claim 
\ref{mt1} of Theorem~\ref{mt}, 
it suffices to prove that one of $X_\al$ is cofinal in $X$. 
 
{\it  Suppose otherwise\/}. 
Then the sets% 
$D_\al=\ens{z\in D}{\sus x\in X_\al\,(z\le x )}$  
contain $\ali$ different sets and form a ROD sequence.

We claim that 
{\it all sets\/ $D_\al$ belong to the same class\/ 
$\fs0\rho$ as the given Borel order\/ $\le$}. 
This will contradict to 
the following lemma by Stern \cite{stern}, and therefore 
complete the proof of item \ref{mt2} of Theorem~\ref{mt}.

\ble 
[in the Solovay model]
\lam{leb}
If\/ $\rho<\omi$ then there is no ROD\/ $\omi$-sequence of 
pairwise different sets\/ $X\sq\bn$ in the class\/ $\fs0\rho$.
\qed
\ele 

To prove the claim, let 
$x_0\le_{\text{lex}}x_1\le_{\text{lex}}x_2\le_{\text{lex}}\ldots$ 
be an arbitrary cofinal sequence in $X_\al$, countable by the 
above. 
Then $D_\al=\ens{z\in D}{\sus n\,(z\le x_n)}$ 
is $\fs0\rho$ by obvious reasons.%

\qeDD{Theorem \ref{mt}\ref{mt2}}

\subsection{$\fs12$ linear suborders of Borel PQOs}

Here we prove Theorem \ref{mt}\ref{mt3}. 
Assume that $\le$ is a Borel PQO on a Borel set $D\sq\bn$,  
$X\sq D$ is a $\fs12$ set, 
and\/ ${\le}\res X$ is a {\it  linear\/} quasi-order.  
We also assume that $\omi^{\rL[r]}<\omi$ for every real $r$. 

Prove that the ordering $\ang{X;{\le}}$ is countably cofinal.

Pick a real $r$ such that $X$ is $\is12(r)$ and 
$\le$ is $\id11(r)$.
To prepare for an absluteness argument, fix canonical formulas,
$$
\bay{rcl}
\vpi(\cdot,\cdot)&\text{of type}& \is12\,,\\[1ex]
\sg(\cdot,\cdot,\cdot)&\text{of type}& \is11\,,\\[1ex]
\pi(\cdot,\cdot,\cdot)&\text{of type}& \ip11\,, 
\eay
$$
which define $\le$ and $X$ in the set universe $\rV$, 
so that it is true
in $\rV$ that
$$
x\le y\leqv \sg(r,x,y)\leqv \pi(r,x,y)
\quad\text{and}\quad
x\in X\leqv\vpi(r,x)\,.
$$
for all $x,y\in\bn.$
We let
$X_{\vpi}=\ens{x\in\bn}{\vpi(r,x)}$
and
$$
x\le_{\sg\pi} y
\quad\eqv\quad \sg(r,x,y)
\quad\eqv\quad \pi(r,x,y)
$$
so that $X_{\vpi}=X$ and $\le_{\sg\pi}$ is $\le$ in $\rV$, but
$X_{\vpi}$ and $\le_{\sg\pi}$ can be defined 
in any transitive universe  containing all ordinals 
(to preserve the equivalence of formulas $\sg$ and $\pi$).
In particular, $X_{\vpi}=X$ and $\le_{\sg\pi}$ is $\le$ 
in the background universe $\rV$.

Let $\wo$ be the canonical $\ip11$ set of codes of (countable) 
ordinals, and for $w\in\wo$ let $|w|<\omi$ be the ordinal 
coded by $w$.

Let $X_\vpi=\bigcup_{\al<\omi}X_\vpi(\al)$ be a canonical 
representation of $X_\vpi$ as an increasing union of $\fs11$ sets.
Thus to define $X_\vpi(\al)$ we fix a $\ip11(r)$ set 
$P\sq(\bn){}^2$ such that $X=\ens{x}{\sus y\,P(x,y)}$, 
fix a canonical $\ip11(r)$ norm 
$f:P\to\omi$, and let 
$$
P_\al=\ens{\ang{x,y}}{f(x,y)<\al}
\quad\text{and}\quad
X_\vpi(\al)=\ens{x}{\sus y\,(\ang{x,y}\in P_\al)}\,.
$$

In our assumptions, the ordinal $\Omega=\omi$ is inaccessible 
in $\rL[r]$. 
Let $\cP=\text{Coll}({{<}\,\Om},\om)\in \rL[r]$ be the corresponding 
Levy collapse forcing.  
Consider a \ddd\cP generic extension 
$\rV[G]$ of the universe.  
Then $\rL[r][G]$ is a Solovay-model generic extension 
of $\rL[r]$.  
The plan is to compare the models $\rV$ and $\rL[r][G]$. 
Note that $\rL[r]$ is their common part, 
$\rV[G]$ is their common extension, and the 
three models have the same cardinal 
$\omi^\rV=\omi^{\rL[r][G]}=\omi^{\rV[G]}=\Om>\omi^{\rL[r]}$. 

By Theorem~\ref{mt}\ref{mt2}, it holds in 
$\rL[r][G]$ that the ordering 
$\stk{X_{\vpi}}{\le_{\sg\pi}}$ is countably cofinal, 
hence there is an ordinal $\al<\Om=\omi^{\rL[r][G]}$ 
such that the sentence 
\ben
\fenu
\itla{f1}
the subset $X_\vpi(\al)$ is 
\ddd{\le_{\sg\pi}}cofinal in the whole set $X_\vpi$ 
\een
is true in $\rL[r][G]$.
However \ref{f1} can be expressed by a $\ip12$ formula with 
$r$ and an arbitrary code $w\in\wo\cap\rL[r][G]$ such that 
$|w|=\al$ --- as the only parameters. 
It follows, by the Shoenfield absoluteness, that \ref{f1} 
is true in $\rV[G]$ as well.

And then, by exactly the same absoluteness argument, 
\ref{f1} is true in the set universe $\rV$, too.
In other words, it is true in $\rV$ that $X_\vpi(\al)$, 
a $\fs11$ set, is cofinal in the whole set $X= X_\vpi$. 
But $X_\vpi(\al)$ is countably cofinal by 
Theorem \ref{mt}\ref{mt1}.\vom

\qeDD{Theorem \ref{mt}\ref{mt3}}

\vyk{
and by absoluteness arguments based on Lemma~\ref{lec}  
we conclude that $X$ is 
ctbly cofinal in the former, too.

\ble 
%[in the Solovay model]
\lam{lec}
Let\/ $\le$ be a Borel PQO on a Borel set\/ $D$, and\/ 
$X\sq D$ a\/ $\fs12$ set. 
Let\/ $\rV^+$ be an extension of the given universe\/ $\rV$ 
with the same\/ $\omi$, and\/ $X^+$ the according extension 
of\/ $X$ in\/ $\rV^+$. 
Then\/ $X$ is cofinal in\/ $X^+$.
\ele 
\bpf[lemma]
A canonical representation of $X=\bigcup_{\al<\omi}X_\al$ as 
a union of $\fs11$ sets reduces the task to the case when $A$ 
is $\fs11$. 
By Theorem~\ref{mt}\ref{mt1} there is, in $V$, a cofinal 
sequence $\vec x=\sis{x_n}{n\in\om}$ in $X$. 
By the Shoenfield absoluteness, it is still true in $\rV^+$  
that $\vec x$ is cofinal in $X^+$.
\epF{lemma}
}

\vyk{

\subsection{Digression: Hausdorff's pantachies}

\bdf 
[Hausdorff 1907, 1909]
\lam{dpant}
A {\it pantachy} is any maximal totally ordered subset $L$ of 
a given partially ordered set $P$, 
\ {\eg,
$P=\stk{\dR^\om}{\le}$}, where, for $x,y\in\dR^\om$, 
$$
x\le y \quad\text{iff}\quad 
x(n)\le y(n) \text{\ \ for all but finite }n, 
$$
the {\it  eventual domination} ordering.   
\edf

\bre
Any pantachy in $P=\stk{\dR^\om}{\le}$ is not 
countably cofinal.
\ere

\bte
[Hausdorff 1909]
\lam{haus1}
There is a pantachy in\/ $\stk{\dR^\om}{\le}$ 
\ {with an\/ $(\omi,\omi)$-gap}.
\ete 

\bprl
[Hausdorff 1907]
Is there a pantachy (in $\stk{\dR^\om}{\le}$), 
\ {containing {\bf no} $(\omi,\omi)$-gaps}\,?
\eprl

The problem is {\bf still open}, and, it looks like it is 
the oldest concrete open problem in set theory.

G\"odel and Solovay discussed almost the same problem in 1970s. 
%(as a matter of fact with no reference to Hausdorff). 

{The problem of effective existence of pantachies}

\bprl
[Hausdorff 1907]
\lam{hp2}
Is there an 
\ {individual, 
effectively defined} example of a pantachy\,?
\eprl

\bsol
%[K \& Lyubetsky 2012]
By Theorem \ref{teb}: 
{\bf generally speaking, in the negative}, 
whenever $P$ is a {Borel}  
partial order, in which
{every countable subset has a strict upper bound}.
\esol

This result, 
by no means surprising, 
is nevertheless based on such a nontrivial argument as 
Stern's theorem. 
But no algebraic structure on $P$ is assumed.

Hausdorff's early papers
}

\end{document}